

  \def \TRUE {Y}
  \def \FALSE {N}

  \def \ifundef #1{\expandafter \ifx \csname #1\endcsname \relax }    



  \def\arxiv{\TRUE}


  \input amssym
  \input miniltx
  \if \arxiv \FALSE
    \input pictex
    \input graphicx.sty
    \input color.sty
    \fi


  \font \bbfive = bbm5
  \font \bbeight = bbm8
  \font \bbten = bbm10
  \font \rs = rsfs10 \font \rssmall = rsfs10 scaled 833  
  \font \eightbf = cmbx8
  \font \eighti = cmmi8 \skewchar \eighti = '177
  \font \fouri = cmmi5 scaled 800 
  \font \eightit = cmti8
  \font \eightrm = cmr8
  \font \eightsl = cmsl8
  \font \eightsy = cmsy8 \skewchar \eightsy = '60
  \font \eighttt = cmtt8 \hyphenchar \eighttt = -1

  \font \sixi = cmmi6 \skewchar \sixi = '177
  \font \sixrm = cmr6
  \font \sixsy = cmsy6 \skewchar \sixsy = '60
  \font \tensc = cmcsc10
  \font \eightsc = cmcsc8
  
  \scriptfont \bffam = \bbeight
  \scriptscriptfont \bffam = \bbfive
  \textfont \bffam = \bbten

  \newskip \ttglue

  \def \eightpoint {\def \rm {\fam 0 \eightrm }\relax
  \textfont 0 = \eightrm \scriptfont 0 = \sixrm \scriptscriptfont 0 = \fiverm
  \textfont 1 = \eighti \scriptfont 1 = \sixi \scriptscriptfont 1 = \fouri
  \textfont 2 = \eightsy \scriptfont 2 = \sixsy \scriptscriptfont 2 = \fivesy
  \textfont 3 = \tenex \scriptfont 3 = \tenex \scriptscriptfont 3 = \tenex
  \def \it {\fam \itfam \eightit }\relax
  \textfont \itfam = \eightit
  \def \sl {\fam \slfam \eightsl }\relax
  \textfont \slfam = \eightsl
  \def \bf {\fam \bffam \eightbf }\relax
  \textfont \bffam = \bbeight \scriptfont \bffam = \bbfive \scriptscriptfont \bffam = \bbfive
  \def \tt {\fam \ttfam \eighttt }\relax
  \textfont \ttfam = \eighttt
  \tt \ttglue = .5em plus.25em minus.15em
  \normalbaselineskip = 9pt
  \def \MF {{\manual opqr}\-{\manual stuq}}\relax
  \let \sc = \sixrm
  \let \big = \eightbig
  \let \rs = \rssmall
  \setbox \strutbox = \hbox {\vrule height7pt depth2pt width0pt}\relax
  \normalbaselines \rm }

  \def \setfont #1{\font \auxfont =#1 \auxfont }
  \def \withfont #1#2{{\setfont {#1}#2}}


  \def \text #1{{\mathchoice {\hbox {\rm #1}} {\hbox {\rm #1}} {\hbox {\eightrm #1}} {\hbox {\sixrm #1}}}}
  \def \varbox #1{\setbox 0\hbox {$#1$}\setbox 1\hbox {$I$}{\ifdim \ht 0< \ht 1 \scriptstyle #1 \else \scriptscriptstyle #1 \fi }}

  \def \rsbox #1{{\mathchoice {\hbox {\rs #1}} {\hbox {\rs #1}} {\hbox {\rssmall #1}} {\hbox {\rssmall #1}}}}
  \def \mathscr #1{\rsbox {#1}}


  \newcount \secno \secno = 0
  \newcount \stno \stno = 0
  \newcount \eqcntr \eqcntr = 0

  \ifundef {showlabel} \global \def \showlabel {\FALSE } \fi
  \ifundef {auxwrite}  \global \def \auxwrite  {\TRUE }  \fi
  \ifundef {auxread}   \global \def \auxread   {\TRUE }  \fi

  \def \undefrule {\kern 2pt \vrule width 2pt height 5pt depth 0pt \kern 2pt}
  \def \UndefLabels {}
  \def \possundef #1{\ifundef {#1}\undefrule {\eighttt #1}\undefrule
    \global \edef \UndefLabels {\UndefLabels #1\par }
  \else \csname #1\endcsname \fi }

  \def \newcommand #1#2{\global \edef #1{#2}} 

  \def \define #1#2{\global \expandafter \edef \csname #1\endcsname {#2}}
  \long \def \error #1{\medskip \noindent {\bf ******* #1}}
  \def \fatal #1{\error {#1\par Exiting...}\end }

  \def \advseqnumbering {\global \advance \stno by 1 \global \eqcntr =0}

  \def \current {\ifnum \secno = 0 \number \stno \else \number \secno \ifnum \stno = 0 \else .\number \stno \fi \fi }


  \def \rem #1{\vadjust {\vbox to 0pt{\vss \hfill \raise 3.5pt \hbox to 0pt{ #1\hss }}}}
  \font \tiny = cmr6 scaled 800
  \def \deflabel #1#2{\relax
    \if \TRUE \showlabel \rem {\tiny #1}\fi
    \ifundef {#1PrimarilyDefined}\relax
      \define {#1}{#2}\relax
      \define {#1PrimarilyDefined}{#2}\relax
      \if \TRUE \auxwrite \immediate \write 1 {\string \newlabe l {#1}{#2}}\fi
    \else
      \edef \old {\csname #1\endcsname }\relax
      \edef \new {#2}\relax
      \ifx \old \new \else \fatal {Duplicate definition for label ``{\tt #1}'', already defined as ``{\tt \old }''.}\fi
      \fi }

  \def \label #1 {\deflabel {#1}{\current }}

  \def \lbldeq #1 $$#2$${\if \InsideProof \FALSE \advseqnumbering
    \edef \lbl {\current }\else
    \global \advance \eqcntr by 1
    \edef \lbl {\current .\number \eqcntr }\fi
    \deflabel {#1}{\lbl }
    $$
    #2
    \eqno {(\lbl )}
    $$}

  \def \split #1.#2.#3.#4;{\global \def \parone {#1}\global \def \partwo {#2}\global \def \parthree {#3}\global \def
\parfour {#4}}
  \def \NA {NA}
  \def \ref #1{\split #1.NA.NA.NA;(\possundef {\parone }\ifx \partwo \NA \else .\partwo \fi )}


  \newcount \bibno \bibno = 0

  \def \Bibitem #1 #2; #3; #4 \par {\smallbreak
    \global \advance \bibno by 1
    \item {[\possundef {#1}]} #2, {``#3''}, #4.\par
    \ifundef {#1PrimarilyDefined}\else
      \fatal {Duplicate definition for bibliography item ``{\tt #1}'',
      already defined in ``{\tt [\csname #1\endcsname ]}''.}
      \fi \ifundef {#1}\else \edef \prevNum {\csname #1\endcsname } \ifnum \bibno =\prevNum \else \error {Mismatch
        bibliography item ``{\tt #1}'', defined earlier (in aux file ?) as ``{\tt \prevNum }'' but should be ``{\tt
        \number \bibno }''.  Running again should fix this.}  \fi \fi
    \define {#1PrimarilyDefined}{#2}\relax
    \if \TRUE \auxwrite \immediate \write 1 {\string \newbi b {#1}{\number \bibno }}\fi }

  \def \jrn #1, #2 (#3), #4-#5;{{\sl #1}, {\bf #2} (#3), #4--#5}
  \def \Article #1 #2; #3; #4 \par {\Bibitem #1 #2; #3; \jrn #4; \par }

  \def \references {\begingroup \bigbreak \eightpoint \centerline {\tensc References} \nobreak \medskip \frenchspacing }


  \catcode `\@ =11
  \def \citetrk #1{{\bf \possundef {#1}}} 
  \def \c@ite #1{{\rm [\citetrk {#1}]}}
  \def \sc@ite [#1]#2{{\rm [\citetrk {#2}\hskip 0.7pt:\hskip 2pt #1]}}
  \def \du@lcite {\if \pe@k [\expandafter \sc@ite \else \expandafter \c@ite \fi }
  \def \cite {\futurelet \pe@k \du@lcite }
  \catcode `\@ =12


  \def \Headlines #1#2{\nopagenumbers
    \headline {\ifnum \pageno = 1 \hfil
    \else \ifodd \pageno \tensc \hfil \lcase {#1} \hfil \folio
    \else \tensc \folio \hfil \lcase {#2} \hfil
    \fi \fi }}

  \font \titlefont = cmbx12
  \long \def \title #1{\begingroup
    \titlefont
    \parindent = 0pt
    \baselineskip = 16pt
    \leftskip = 35truept plus 1fill
    \rightskip = \leftskip
    #1\par \endgroup }

  \long \def \Quote #1\endQuote {\begingroup \leftskip 35truept \rightskip 35truept \parindent 17pt \eightpoint #1\par \endgroup }
  \long \def \Abstract #1\endAbstract {\vskip 1cm \Quote \noindent #1\endQuote }

  \def \Note #1{\footnote {}{\eightpoint #1}}
  \def \Date #1 {\Note {\it Date: #1.}}

  \newcount \auxone \newcount \auxtwo \newcount \auxthree
  \def \currenttime {\auxone =\time \auxtwo =\time \divide \auxone by 60 \auxthree =\auxone \multiply \auxthree by 60
\advance
    \auxtwo by -\auxthree \ifnum \auxone <10 0\fi \number \auxone :\ifnum \auxtwo <10 0\fi \number \auxtwo }
  \def \today {\ifcase \month \or January\or February\or March\or April\or May\or June\or July\or August\or September\or
    October\or November\or December\fi { }\number \day , \number \year }
  
  \def \hojeExtenso {\number \day \ de \ifcase \month \or janeiro\or fevereiro\or mar\c co\or abril\or maio\or junho\or
julho\or
    agosto\or setembro\or outubro\or novembro\or decembro\fi \ de \number \year }

  \def \part #1#2{\vfill \eject \null \vskip 0.3\vsize
    \withfont {cmbx10 scaled 1440}{\centerline {PART #1} \vskip 1.5cm \centerline {#2}} \vfill \eject }


  \def \Case #1:{\medskip \noindent {\tensc Case #1:}}

  \def \fix {\smallskip \noindent $\blacktriangleright $\kern 12pt}
  \def \newpage {\vfill \eject }

  \def \lcase #1{\edef \auxvar {\lowercase {#1}}\auxvar }

  \def \emph #1{{\it #1}\/}

  \def \state #1 #2\par {\begingroup \def \InsideProof {\TRUE} \medbreak \noindent \advseqnumbering {\bf \current .\enspace
#1.\enspace \sl #2\par }\medbreak \endgroup }

  \def \lstate #1 #2 #3\par {\state {#2} \label {#1} #3\par }

  \def \definition #1\par {\state Definition \rm #1\par }




  \def \Proof {\if \InsideProof \FALSE \else \fatal {Opening Proof block before closing the last one} \fi
    \global \def \InsideProof {\TRUE}
    \medbreak \noindent {\it Proof.\enspace }}

  \def \endProof {\if \InsideProof \TRUE \hfill $\endproofmarker $ \looseness = -1 \fi
    \medbreak
    \if \InsideProof \FALSE \fatal {Closing Proof block before opening it} \fi
    \global \def \InsideProof {\FALSE}}

  \def \closeProof {\eqno \endproofmarker
    \global \def \InsideProof {E}}

  \def \InsideProof {\FALSE}


  \def \quebra #1{#1 $$$$ #1}
  \def \explica #1#2{\mathrel {\buildrel \hbox {\sixrm #1} \over #2}}
  \def \explain #1#2{\explica {\ref {#1}}{#2}} 
  
  \def \=#1{\explain {#1}{=}}

  \newcount \fnctr \fnctr = 0
  \def \fn #1{\global \advance \fnctr by 1
    \edef \footnumb {$^{\number \fnctr }$}\relax
    \footnote {\footnumb }{\eightpoint #1\par \vskip -10pt}}

  \def \text #1{\mathchoice {\hbox {#1}} {\hbox {#1}} {\hbox {\eightrm #1}} {\hbox {\sixrm #1}}}


  \def \item #1{\par \noindent \kern 1.1truecm\hangindent 1.1truecm \llap {#1\enspace }\ignorespaces }
  
  \def \Item #1{\smallskip \item {{\rm #1}}}

  \newcount \zitemno \zitemno = 0
  \def \izitem {\global \zitemno = 0}

  \def \zitemplus {\global \advance \zitemno by 1 \relax }
  \def \rzitem {\romannumeral \zitemno }
  \def \rzitemplus {\zitemplus \rzitem }
  \def \zitem {\Item {{\rm (\rzitemplus )}}}

  \def \lzitem #1 {\Item {{\rm (\rzitemplus )}} {\deflabel {#1}{\current .\rzitem }{\def \showlabel {\FALSE }\deflabel
{Local#1}{\rzitem }}}}

  \newcount \nitemno \nitemno = 0
  
  \def \nitem {\global \advance \nitemno by 1 \Item {{\rm (\number \nitemno )}}}

  \newcount \aitemno \aitemno = -1
  \def \boxlet #1{\hbox to 6.5pt{\hfill #1\hfill }}
  
  \def \aitemconv {\ifcase \aitemno a\or b\or c\or d\or e\or f\or g\or h\or i\or j\or k\or l\or m\or n\or o\or p\or q\or
    r\or s\or t\or u\or v\or w\or x\or y\or z\else zzz\fi }
  \def \aitem {\global \advance \aitemno by 1\Item {(\boxlet \aitemconv )}}

  \def \bitem {\Item {$\bullet $}}


  \def \deflabeloc #1#2{\deflabel {#1}{\current .#2}{\edef \showlabel {\FALSE }\deflabel {Local#1}{#2}}}
  \def \lbldzitem #1 {\zitem \deflabeloc {#1}{\rzitem }}
  \def \lbldaitem #1 {\aitem \deflabeloc {#1}{\aitemconv }}
  \def \aitemmark #1 {\deflabel {#1}{\aitemconv }}

  \def \zitemmark #1 {\deflabel {#1}{\current .\rzitem }{\def \showlabel {\FALSE }\deflabel {Local#1}{\rzitem }}}


  \def \mathbb #1{{\bf #1}}
  \def \text #1{\hbox {#1}}
  \def \frac #1#2{{#1\over #2}}

  \def \<{\left \langle \vrule width 0pt depth 0pt height 8pt }
  \def \>{\right \rangle }
  \def \({\big (}
  \def \){\big )}
  
  \def \and {\mathchoice {\hbox {\quad and \quad }} {\hbox { and }} {\hbox { and }} {\hbox { and }}}

  \def \imply {\mathrel {\Rightarrow }}
  \def \IFF {\kern 7pt\Leftrightarrow \kern 7pt}
  \def \IMPLY {\kern 7pt \Rightarrow \kern 7pt}
  
  \def \endproofmarker {\square }
  \def \"#1{{\it #1}\/} 
  
  \def \*{\otimes }
  \def \caldef #1{\global \expandafter \edef \csname #1\endcsname {{\cal #1}}}
  \def \mathcal #1{{\cal #1}}
  \def \bfdef #1{\global \expandafter \edef \csname #1\endcsname {{\bf #1}}}
  \bfdef N \bfdef Z \bfdef C \bfdef R
  \def \exists {\mathchar "0239\kern 1pt }
  \def \labelarrow #1{\setbox 0\hbox {\ \ $#1$\ \ }\ {\buildrel \textstyle #1 \over {\hbox to \wd 0 {\rightarrowfill }}}\ }
  \def \subProof #1{\medskip \noindent #1\enspace }
  \def \itmProof (#1) {\subProof {(#1)}}
  \def \itemImply #1#2{\subProof {#1$\Rightarrow $#2}}
  \def \itmImply (#1) > (#2) {\itemImply {(#1)}{(#2)}}


  \if \TRUE \auxread
    %
    \IfFileExists {\jobname .aux}{\input \jobname .aux}{\null } \fi
  \if \TRUE \auxwrite \immediate \openout 1 \jobname .aux \fi


  \def \close {\ifx \empty \UndefLabels \else
    \message {*** There were undefined labels ***} \medskip \noindent
    ****************** \ Undefined Labels: \tt \par \UndefLabels \fi
    \if \TRUE \auxwrite \closeout 1 \fi
    \par \vfill \supereject \end }


  \def \Caixa #1{\setbox 1=\hbox {$#1$\kern 1pt}\global \edef \tamcaixa {\the \wd 1}\box 1}
  \def \caixa #1{\hbox to \tamcaixa {$#1$\hfil }}

  \def \medcup {\mathop {\mathchoice {\raise 1pt \hbox {$\mathchar "1353$}}{\mathchar "1353}{\mathchar "1353}{\mathchar
"1353}}}
  \def \medcap {\mathop {\mathchoice {\raise 1pt \hbox {$\mathchar "1354$}}{\mathchar "1354}{\mathchar "1354}{\mathchar
"1354}}}


  \def \cl #1 #2 #3 {#1, & \hbox {#2 } #3\hfill \crr }

  \def \paper #1#2#3{
    
    \hsize #1truemm   \advance \hsize by -#3truemm  \advance \hsize by -#3truemm
    \vsize #2truemm   \advance \vsize by -#3truemm  \advance \vsize by -#3truemm
    \hoffset =-1truein \advance \hoffset by #3truemm
    \voffset =-1truein \advance \voffset by #3truemm
    }

  %
  %

  \def \section #1 \par {\global \advance \secno by 1 \stno = 0
    \goodbreak \bigbreak
    \noindent {\bf \number \secno .\enspace #1.}
    \nobreak \medskip \noindent }

  \newcount \loopVariable
  \def \massagem #1{
    \loopVariable = 200 \loop \message {*} \advance \loopVariable by-1 \ifnum \loopVariable >0 \repeat
    \message {#1}
    \loopVariable = 200 \loop \message {*} \advance \loopVariable by-1 \ifnum \loopVariable >0 \repeat
    }

  \def \startsection #1 \par
    {\goodbreak
    \bigbreak
    \if \arxiv \FALSE
      \ifnum \currentgrouplevel > 0
        \vfill \eject \null \vfill
        \massagem {*** DO NOT START SECTION INSIDE A GROUP, SUCH AS BEFORE CLOSING LAST ONE ***}
        *** DO NOT START SECTION INSIDE A GROUP, SUCH AS BEFORE CLOSING LAST ONE ***
        \vfill \eject
        \fi
      \fi
    \begingroup
    \global \edef \secname {#1}\relax
    \global \advance \secno by 1
    \stno = 0
    \AddToTableOfContents {\number \secno .}{#1}}

  \def \sectiontitle \par
    {\noindent {\bf \number \secno .\enspace \secname .}
    \nobreak \medskip
    \noindent }

  \def \endsection {\endgroup }

  %
  %

  \def \ToCf {\TRUE }
  \if \ToCf \TRUE \def \AddToTableOfContents #1#2{{\let \the =0\edef \a {\write 2{\string \inde x #1; #2; \the \count 0;}}\a }}
            \else \def \AddToTableOfContents #1#2{}\fi
  \def \condinput #1 {\IfFileExists {#1}{\input #1}{*** MISSING FILE #1.  RUNNING AGAIN MIGHT FIX THIS.}}
  \def \index #1; #2; #3;{\line {\hbox to 20pt{\hfill #1} #2 \ \dotfill \ \ #3}}
  \def \tableOfContents {\centerline {CONTENTS} \bigskip \IfFileExists {c}{\input c}{*** MISSING FILE c.  RUNNING AGAIN MIGHT FIX THIS.}ontents.aux \if \ToCf \TRUE \openout 2 contents.aux \fi }

  %
  %

  \def \dotedline #1 ... #2{\smallskip \line {\qquad #1\quad \dotfill \quad #2\qquad }}
  \def \PrintSymbol #1; #2; #3;{\dotedline \hbox to 1.5cm{\hfill $#1$\hfill }\ #2 ... {#3}}


  \def \PrintConcept #1; #2;{\dotedline #1 ... {#2}}

  \def \ConceptBank {}

  \def \newConcept #1#2{\emph {#1}\def \tempVar {#2}\ifx \empty \tempVar \def \tempVar {#1}\else \def \tempVar {#2}\fi
    \global \edef \ConceptBank {\ConceptBank \PrintConcept \tempVar ; \number \secno .\number \stno ;}}

  \def \PART #1{\newpage \setbox 0\hbox {\bf PART #1}
    \dimen 0 \wd 0
    \advance \dimen 0 by 20pt
    \centerline {\hbox to \dimen 0{\hrulefill }}
    \medskip \centerline {\box 0}
    \centerline {\hbox to \dimen 0{\hrulefill }}
    \vskip 2cm}

\def \imply {\dimen 0=3pt\kern \dimen 0\Rightarrow \kern \dimen 0}
\def \iff {\dimen 0=3pt\kern \dimen 0\Leftrightarrow \kern \dimen 0}

\def \paren #1{\setbox 0\hbox {$#1$}\ifdim \ht 0 < 7.17pt (#1) \else \big (#1\big )\fi }

\def \axi /{{\eightsc (inv)}}

\def \invr {{\text {\eightrm inv}}}
\def \I {{\mathscr I}}
\def \R {{\mathscr R}}
\def \Iinv {\I_\invr}
\def \Rinv {\R_\invr}
\def \Rprimeinv {\R\,'_\invr}

\def \cl #1{\overline {#1}}

\def \lanFun {\text {Ann}^{\scriptscriptstyle L}}
\def \ranFun {\text {Ann}^{\scriptscriptstyle R}}

\def \lan #1#2{\lanFun _{\scriptscriptstyle #1}\paren {#2}}
\def \ran #1#2{\ranFun _{\scriptscriptstyle #1}\paren {#2}}
\def \an #1#2{\text {Ann}_{\scriptscriptstyle #1}\paren {#2}}

\def \lanA #1{\lan {}{#1}}
\def \ranA #1{\ran {}{#1}}
\def \anA #1{\an {}{#1}}

\def \idealGen #1#2{\langle #2\kern 1pt\rangle _{\kern -1pt\scriptscriptstyle #1}}
\def \idealGenA #1{\idealGen {}#1}

\title {\bf REGULAR IDEALS UNDER THE IDEAL\par INTERSECTION PROPERTY}

\bigskip
\centerline {\tensc
  R. Exel\footnote {$^{\ast }$}{\eightrm Departamento de Matem\'atica -- Universidade Federal de Santa Catarina.}
  }

\Note {Partially supported by CNPq.}

\Abstract
  The goal of this short note is to prove that when $A$ is a closed *-subalgebra of a C*-algebra $B$ satisfying the
ideal intersection property plus a mild axiom \axi /, then the map $J\mapsto J\cap A$  establishes an isomorphism from the
boolean algebra of all regular ideals of $B$ to the boolean algebra of all regular, invariant ideals of $A$.
 \endAbstract

\startsection Introduction

\def \G {G}

\sectiontitle

Recall that a closed *-subalgebra $A$ of a C*-algebra $B$ is said to satisfy the \emph {ideal intersection property} if
every nonzero ideal of $B$ has a  nonzero intersection with $A$.

This property has a long history, partially motivated by the celebrated uniqueness Theorem of Cuntz and Krieger
\cite [Theorem 2.13]{CK}, originally stated in a different way but which can be directly linked to the ideal intersection
property.  In a more explicit fashion, Renault identified  a situation \cite [Proposition II.4.6]{Renault}  where
the ideals of $B$ and the ideals of  $A$ are even more closely related,
namely the correspondence given by
  \lbldeq RenaultCorresondence
  $$
  J \in \I (B) \mapsto J\cap A \in \Iinv (A),
  $$
  is a bijection,
  where $\I (B)$ denotes the collection of all ideals of $B$, and $\Iinv (A)$ denotes the collection of all
\emph {invariant ideals} of $A$ (see below for the appropriate definitions).

  The conditions under which this result was proven by Renault consist in taking  $B$ to be the reduced C*-algebra of a twisted \'etale groupoid
$\G $, and taking  $A$ to be the abelian subalgebra
  associated to the unit space $\G ^{(0)}$.
  It is also required that $\G $ be an
\emph {essentially principal} groupoid \cite [Definition II.4.3]{Renault}, in the sense that for every closed invariant
subset $F\subseteq \G ^{(0)}$, one has that the set of points in $F$ with trivial isotropy is dense in $F$.

A weaker notion of essential principality has been considered by several authors, Renault included, where the
above requirement is made just for $F=\G ^{(0)}$, rather than for all closed invariant subsets (see \cite [Proposition
3.1]{RenaultCartan}).

This weaker version of essentially principal groupoids is not strong enough to guarantee \ref {RenaultCorresondence} to
be a bijection, but one may still recover the ideal intersection property, a fact I believe was first proved in
\cite [Theorem 4.4]{ExelNHausd} (see also \cite [Theorem 3.1.(b)]{Lots}).  In fact, in the second-countable case, the ideal
intersection property turns out to be equivalent to the above weaker form of essential principality \cite [Proposition 5.5.(2)]{BCFS}.

In the  recent paper \cite {RegIdeals} the authors take an intriguing middle of the road point of view, assuming
the ideal intersection property (same as the weaker form of essential principality in the groupoid case,  as already
observed), and aiming
for the conclusion of the  stronger result proved by Renault mentioned above  (which requires the full essential
principality property), with the caveat that only regular ideals are considered, namely those coinciding with their
second anihilator.  To be precise, \cite [Theorem 3.18]{RegIdeals} states that the correspondence
  \lbldeq PittsCorresondence
  $$
  J \in \R (B) \mapsto J\cap A \in \Rinv (A),
  $$
  is a bijection,
  where $\R (B)$ denotes the collection of all regular ideals of $B$, and $\Rinv (A)$ denotes the collection of all
regular, invariant ideals of $A$.  Besides assuming the ideal intersection property,
this result requires that $A$ be a regular subalgebra of $B$ (in the sense that $B$ is spanned by the normalizers of
$A$), and that $B$ admits an invariant faithful conditional expectation  onto $A$ (see \cite {RegIdeals}  for more
details).

The purpose of the present work is to generalize the above result
  by proving that \ref {PittsCorresondence} is a
bijection under much weaker hypothesis.  Besides the obvious assumption that $A$ satisfy the ideal intersection
property, we need only require axiom \axi /, a weak form of regularity for subalgebras  (see
\ref {DefineAxi} below for the precise definition)
and, more importantly, no assumptions are made on the existence of conditional
expectations whatsoever.
Our main result in this direction is Theorem  \ref {MainThm}.

We then specialize to the context of regular subalgebras in Corollary \ref {MainCorol}, where the only hurdle to be overcome is to check that
axiom \axi / follows from regularity.

I'd like to express my thanks to David Pitts for many interesting conversations regarding regular ideals during his
recent visit to Florian\'opolis, when he pointed me to the interesting preprint
 \cite {RegIdeals} by himself and collaborators.

\endsection

\startsection Anihilators

\sectiontitle

In this section we shall fix our notation regarding anihilators of subsets of a C*-algebra and we will  also introduce a notion
of \emph {invariance} which is important for the formulation of our main result and which we have not found in the literature.

\fix Throughout this section $A $ will  denote a fixed C*-algebra.

\definition
Given any subset $S\subseteq A $,
  \izitem
  \zitem the \emph {left-annihilator of $S$ in $A $} is defined by
  $$
  \lan A S = \{x\in A : xs=0, \ \forall s\in S\}.
  $$
  \zitem the \emph {right-annihilator of $S$ in $A $} is defined by
  $$
  \ran A S = \{x\in A : sx=0, \ \forall s\in S\}.
  $$
  \smallskip \noindent In case there is no question as to which ambient algebra we are referring, we will drop the
subscript and write simply $\lanA S$ and $\ranA S$.

\state Proposition \label ElemStuff
  Given a subset $S\subseteq A $, one has that:
  \izitem
  \zitem $\lanA S$ and $\ranA S$ are norm-closed,
  \zitem $\lanA S$ is always a left ideal,
  \zitem if $S$ is a left ideal then $\lanA S$ is a right (and hence two-sided) ideal,
  \zitem $\ranA S$ is a always a right ideal,
  \zitem if $S$ is a right ideal then $\ranA S$ is a left  (and hence two-sided) ideal.

\Proof Left for the reader.
\endProof

Recall that,
given any subset $S\subseteq A $, it is customary to denote its closed linear span by $[S]$.  In symbols:
  $$
  [S] = \overline {\text {span}(S)}.
  $$
  If $T$ is another subset of $A $,
the following notation will also be used
  $$
  ST=\{st: s\in S, \ t\in T\}.
  $$
  The closed two-sided ideal generated by $S$ in $A $ will be denoted by $\idealGen A S$, so that
  $$
  \idealGen A S = [A SA ].
  $$
  As before we will write   $\idealGenA S$ when there is no chance of confusion.

\definition We will say that a subset $S\subseteq A $ is \emph {$A $-invariant} provided $[SA ]=[A S]$.

See Section \ref {RegularSection}, below, for a discussion of the above concept in the context of regular inclusions.

  Notice that $[SA ]$ is always a right ideal and that $[A S]$ is always a left ideal.  Thus, when  $S$ is $A $-invariant, one
has that $[SA ]$ is a closed two-sided ideal, clearly the smaller such ideal containing
$S$,  whence
  \lbldeq RightLeftTwo
  $$
  [SA ]=\idealGenA S = [A S].
  $$
  In addition,  when  $S$ is a closed two-sided ideal, it is easy to see that $[SA ]=S=[A S]$, so  $S$ is $A $-invariant.

\state Proposition \label AnnAndIdeals
  If  $S$ is $A $-invariant  then
  $$
  \lanA S = \ranA S = \lanA {\idealGenA S} = \ranA {\idealGenA S}.
  $$

\Proof
  Suppose first that $S$ is a closed two-sided ideal,  therefore also self-adjoint.  Then, given $x\in \lanA S$, and $s\in S$, we have that $x^*s^* \in S$, whence
  $$
  (sx)  (sx)^*  = sxx^*s^* \in sxS= \{0\}.
  $$
  This shows that $sx=0$, and hence that $x\in \ranA S$, proving that  $\lanA S \subseteq \ranA S$, while the reverse inclusion follows similarly.
  In the general case
  observe that
  $$
  \lanA S = \lanA {SA } = \lanA {[SA ]} \={RightLeftTwo} \lanA {\idealGenA S},
  $$
  and similarly
  $$
  \ranA S = \ranA {A S} = \ranA {[A S]} \={RightLeftTwo} \ranA {\idealGenA S}.
  $$
  The final  conclusion then follows from the first part of the proof.
\endProof

\definition \label definePerp
  If $S$ is $A $-invariant (e.g.~a closed two-sided ideal), we let the \emph {annihilator of $S$ in $A $} be defined
by
  $$
  \an A S := \lan A S = \ran A S.
  $$

  We remark that some authors use $S^\perp $ for $\an AS$, but we shall prefer the latter notation due to the fact that we
will soon consider another C*-algebra $B$ containing $A$, and it will therefore be important to specify the ambient algebra
relative to which we are considering the annihilator of a subset $S\subseteq A$.  In this case we feel that $\an AS$ and $\an BS$
are preferable to  $S^{\perp _A}$
and $S^{\perp _B}$.
  Meanwhile, when there is no question as to which ambient algebra we are considering, we will drop the
subscript and write simply $\anA S$.

Notice that,  under the above assumption that $S$ is $A $-invariant,  $\anA S$ is a closed two-sided ideal of $A $.
  By \ref {AnnAndIdeals} we
moreover have that
  \lbldeq ppGenIdeal
  $$
  \anA S= \anA {\idealGenA S}.
  $$

\definition
  An ideal\fn {Unqualified references to \emph {ideals} will henceforth mean closed two-sided
ideals.} $J\trianglelefteq A$ is said to be regular if
  $$
  \anA {\anA J} =J.
  $$

\state Lemma \label LemmaAnnihil
  Let $S$ be an $A$-invariant subset of $A$.  Then:
  \izitem
  \zitem $\anA {\anA {\anA S}} = \anA S$,
  \zitem $\anA S$ is a regular ideal of $A$,
  \zitem If $J_1$ and $J_2$ are regular ideals,  then $J_1\cap J_2$ is regular.

\Proof
Consider the symmetric relation ``$\kern -1pt\perp $\kern -1pt" defined on $A$ by
  $$
  x\perp y \kern 6pt \Leftrightarrow \kern 6pt xy=0=yx.
  $$
  We then claim that
  $$
  S^\perp := \{x\in A: x\perp s, \ \forall s\in S\} =\an AS.
  $$
  Indeed, by \ref {definePerp} we have that
  $$
  \an AS = \lan AS \cap \ran AS  \quebra =
  \{x\in A: xs=0, \ \forall s\in S\} \kern 4pt \cap \kern 4pt \{x\in A: sx=0, \ \forall s\in S\} \quebra =
  \{x\in A: x\perp s=0, \ \forall s\in S\}  = S^\perp ,
  $$
  proving the claim.  Employing \ref {TongueTwisterTwo.ii}, below,  we then have that $S^\perp =((S^\perp )^\perp )^\perp $, which translates into
  (i) and in turn implies (ii).

Regarding (iii),  we have by \ref {TongueTwisterTwo.i}  that
  $$
  J_1\cap J_2 \subseteq \anA {\anA {J_1\cap J_2}}.
  $$
  On the other hand we have that   $J_1\cap J_2 \subseteq J_1$, so
  $$
  \anA {\anA {J_1\cap J_2}}\subseteq \anA {\anA {J_1}} = J_1,
  $$
  and likewise   $\anA {\anA {J_1\cap J_2}} \subseteq J_2$, so
  $$
  \anA {\anA {J_1\cap J_2}}\subseteq J_1\cap J_2.
  $$
  This proves that $J_1\cap J_2$ is regular.
\endProof

\state Lemma \label TongueTwisterTwo
  Let $X$ be any set  and let ``$\perp $\kern -1pt" be any symmetric relation on $X$.  For each subset $S\subseteq X$ define
  $$
  S^\perp = \{x\in X: x\mathrel {\perp } s, \ \forall s\in S\}.
  $$
  Then
  \izitem
  \zitem
  $ S \subseteq (S^\perp )^\perp $, and
  \zitem $S^\perp = ((S^\perp )^\perp )^\perp $.

\Proof Despite sounding a bit like a tongue-twister  it is obvious that every  element in  $S$ is perpendicular  to
anything that is perpendicular to every element in $S$,  hence (i).
  Substituting $S$ for $S^\perp $ in  (i), we get $S^\perp \subseteq ((S^\perp )^\perp )^\perp $.
Next observe that
  $$
  S_1\subseteq S_2 \mathrel {\Rightarrow }  S_2^\perp \subseteq S_1^\perp ,
  $$
  which we may apply to   (i), leading to
  $((S^\perp )^\perp )^\perp \subseteq S^\perp $, and concluding the proof.
  \endProof

\medskip
Observe that \ref {LemmaAnnihil.iii} implies that the collection of all regular ideals of $A$ is a meet-semi-lattice under the
operation
  $$
  J_1 \wedge J_2 := J_1 \cap J_2,
  $$
  in the sense that  $J_1 \cap J_2$ is the largest regular
  ideal\fn {Actually the \emph {largest ideal}.}
  smaller than both $J_1$ and $J_2$.

There is also a smallest regular
  ideal\fn {As opposed to the \emph {smallest ideal},  which is $J_1+J_2$, but which is not always regular.}
  larger than both $J_1$ and $J_2$, namely
  $$
  J_1 \vee J_2 :=
  \anA {\anA {J_1+J_2}},
  $$
  so the collection of all regular ideals is actually a lattice.  If we furthermore define the \emph {negation} of a
regular ideal $J$ by
  $$
  \neg J = \anA J,
  $$
  we get a boolean algebra which may be proved to be complete.

\medskip

As seen above, any ideal $I$ of $A$ is automatically $A$-invariant, but in case $A$ is a closed *-subalgebra of another
C*-algebra $B$,  there is no reason for $I$  to be also $B$-invariant.  In the most
favorable case we have:

\state Proposition  \label HandyPropNew
  Let $I$ be an ideal of $A$ which is also $B$-invariant,  where $B$ is a C*-algebra with $A\subseteq B$.  Then
  $$
  \an AI = \an BI\cap A.
  $$

\Proof
We have that
  $$
  \an A {I} =   \lan A {I} = \lan B {I} \cap A = \an B {I} \cap A.
  \closeProof
  $$
\endProof

\endsection

\startsection Ideal intersection property

\def \standingHyp {the ideal intersection property and axiom \axi /}

\sectiontitle

\fix From now on we will fix a C*-algebra $B$ and a closed *-subalgebra $A\subseteq B$.

\definition We shall  say that $A$ satisfies the \emph {ideal intersection property (relative to $B$)} if, for any nonzero
  ideal
  $J\trianglelefteq B$, one has that $J\cap A\neq \{0\}$.

The following is a key technical result.  It is essentially the only place where we will explicitly invoke the ideal
intersection property.

\state Lemma \label InterJAPerp
  Assume that $A$ satisfies the ideal intersection property, and let $J\trianglelefteq B$ be an ideal.  Then
  $\an B{\idealGen B {J\cap A}} = \an B J$.

\Proof
  Given $b\in \an B{\idealGen B {J\cap A}} $, let
$K=\idealGen B b \cap J$.  We then have that
  $$
  K\cap A=\idealGen B b \cap J \cap A \subseteq
  \an B{\idealGen B {J\cap A}} \cap \idealGen B {J\cap A} =
  \{0\},
  $$
  whence $K=\{0\}$ by the ideal intersection property.  Therefore
  $$
  bJ\subseteq \idealGen B b \cdot J =  \idealGen B b \cap J = K = \{0\},
  $$
  showing that $b\in \lan B J = \an B J $.
This proves that
$\an B{\idealGen B {J\cap A}} \subseteq \an B J $.  The reverse inclusion is clear. \endProof

We now wish to replace the occurrence of $\an B{\idealGen B {J\cap A}}$ in the above Lemma with $\an B{J\cap A}$ but, unless we know that
$J\cap A$ is $B$-invariant, we cannot even make sense of  $\an B{J\cap A}$.

\definition \label DefineAxi
  We shall say that $A$ satisfies \emph {axiom \axi / relative to $B$} when  $J\cap A$
is $B$-invariant  for every ideal $J\trianglelefteq B$.

Please see the following section for more context on this axiom, where it is shown that \axi /  holds for every regular subalgebra.

\state Corollary \label TwoInOne
  Assume that  $A$ satisfies \standingHyp .  Then, for every ideal
$J\trianglelefteq B$, one has that
  \izitem
  \zitem
  $
  \an B{J\cap A} = \an B J,
  $
  and
  \zitem
  $
  \an A{J\cap A} =  \an BJ \cap A.
  $

\Proof
  In order to prove (i),  observe that $J\cap A$ is $B$-invariant by \axi /, so
  $$
  \an B{J\cap A} \={ppGenIdeal} \an B{\idealGen B {J\cap A}} \={InterJAPerp} \an B J.
  $$
  Regarding (ii), we have that
  $$
  \an A{J\cap A} \={HandyPropNew}
  \an B{J\cap A} \cap A \explica {(i)}=
  \an BJ \cap A.
  \closeProof
  $$
  \endProof

The following is our main result.

\state Theorem \label MainThm
  Let $B$ be a C*-algebra and let $A$ be a closed *-subalgebra of $B$ satisfying \standingHyp .
  Also let:
  \bitem $\R(B)$ be the family formed by all regular ideals of $B$, and
  \bitem $\Rinv (A) $ be the family formed by all regular, $B$-invariant ideals of $A$.
  \medskip \noindent Then the correspondence
  $$
  \alpha :J\in \R(B)\ \mapsto \ J\cap A\in \Rinv (A)
  $$
  is a boolean algebra isomorphism, and its inverse is given by
  $$
  \beta :I\in \Rinv (A) \ \mapsto \ \an B{\an B{I}}\in \R(B).
  $$

\Proof We first show that $\alpha $ does indeed map $\R(B)$ into $\Rinv (A) $.  Given a regular ideal $J\trianglelefteq B$, we have that
$J\cap A$ is $B$-invariant by axiom \axi /, so it remains to prove that it is also regular as an ideal of $A$.  For this we
need to show that
  \lbldeq JAReg
  $$
  \an A{\an A{J\cap A}} = J\cap A.
  $$
  Plugging \ref {TwoInOne.ii} into the left-hand-side above  yields
  $$
  \an A{\an A{J\cap A}} =
  \an A{  \an BJ \cap A} \explica {(!)}=   \an B{\an BJ} \cap A = J\cap A,
  $$
  where the passage marked with $(!)$ follows by substituting $\an BJ$ for $J$ in   \ref {TwoInOne.ii}.  This proves \ref {JAReg}.

Observe that, for every $B$-invariant regular ideal
$I\trianglelefteq A$,  one has that $\beta (I)$ is a regular ideal of $B$ by
 \ref {LemmaAnnihil.ii}.  Therefore  $\beta $ maps $\Rinv (A) $ into $\R(B)$.

Given  $I$ in $\Rinv (A) $, we next show that $\alpha \big (\beta (I)\big )=I$.  Because $I$ is regular in $A$,  we have that
  $$
  I =
  \an A {\an AI} \={HandyPropNew}
  \an A {\an BI\cap A} \={TwoInOne.ii}
  \an B{\an BI} \cap A = \alpha (\beta (I)).
  $$
On the other hand,
given a regular ideal $J\trianglelefteq B$,  we have
  $$
  \beta \big (\alpha (J)\big )=
  \an B{\an B{J\cap A}}  \={TwoInOne.i}
  \an B{\an BJ}  =
  J.
  $$
  This shows that $\alpha $ and $\beta $ are each other's inverse.   It is also clear that both $\alpha $ and $\beta $ preserve inclusion of
ideals, so $\alpha $ is a boolean algebra isomorphism, as required.
  \endProof

\medskip As an aside, observe that the map $\beta $ referred to in \ref {MainThm} may be equivalently defined by
  $$
  \beta (J)=\an B{\an A{I}},
  $$
  because
  $$
  \an B{\an A{I}} \={HandyPropNew}
  \an B{\an B{I} \cap A} \={TwoInOne.i}
  \an B{\an B{I}}.
  $$

\endsection

\startsection Regular inclusions

\label RegularSection

\sectiontitle

Recall that axiom \axi / was introduced in  \ref {DefineAxi}  as a solution to the problem as to whether or not $J\cap A$ is
$B$-invariant for every ideal $J\trianglelefteq B$.

Postulating away one's problems might of course sound like a cheap solution, so we will now attempt to dispel this
impression by identifying a very common situation (Proposition \ref {JAIsBInvar.iii}, below) in which the above axiom is
automatically satisfied.

\fix As before, we fix a C*-algebra $B$ and a closed *-subalgebra $A\subseteq B$.

\medskip Recall from \cite [2.3.1]{ExelPitts} that a \emph {normalizer} of $A$ in $B$ is any element $n$ of $B$ such that both $n^*An$ and
$nAn^*$ are contained in $A$.   When
  the normalizers span a dense subspace of $B$,
  and  $A$ contains an approximate unit  for $B$,
  we say that $A$ is a \emph {regular
    subalgebra}\fn {Despite the formal similarity there is no relationship between the concepts of \emph {regular
    subalgebras} and \emph {regular ideals}.  The former was perhaps first introduced in \cite {Dixmier} in the context of von
    Neumann algebras, while the latter is related to the concept of \emph {regular open sets}, namely sets coinciding with
    the interior of its closure.}.

\medskip
Besides the concept of $B$-invariance,  there is another notion of invariance frequently featuring in works dealing with
regular subalgebras, as follows:

\definition A subset $S\subseteq A$ will be called  \emph {normalizer-invariant} when  $nSn^*\subseteq S$, for every normalizer $n$.

The two notions of invariance are closely related to each other, as we now show.

\state Proposition \label JAIsBInvar
  Assume that $A$ is a regular subalgebra of $B$.  Then
  \izitem
  \zitem for every ideal $J\trianglelefteq B$,  one has that $J\cap A$ is normalizer-invariant,
  \zitem every normalizer-invariant ideal of $A$ is $B$-invariant,
  \zitem $A$ satisfies axiom \axi / relative to $B$.

\Proof
  Point (i) is trivial and hence is left for the reader.  As for (ii), let $I$ be a normalizer-invariant ideal of $A$.
Given a normalizer $n$, and given $a$ in $I$, we first claim that
  \lbldeq naInIB
  $$
  na\in [IB].
  $$
  To prove this it is clearly enough to show that, choosing an approximate unit $\{e_i\}_i$ for $I$, one has that
  $$
  na =  \lim _{i\to \infty }e_ina.
  $$
  Let us therefore compute
  $$
  \|e_ina-na\|^2 =
  \|(e_ina-na)(e_ina-na)^*\|
  \quebra =
  \|e_inaa^*n^*e_i  - e_inaa^*n^*-  naa^*n^*e_i  + naa^*n^*\| \quebra \leq
  \|e_i\|\|naa^*n^*e_i  - naa^*n^*\|+ \| naa^*n^*e_i  - naa^*n^*\|,
  $$
  which converges to zero as $i\to \infty $,  because $naa^*n^*\in nIn^*\subseteq I$,  hence proving \ref {naInIB}.
  Since the linear span of the  normalizers is dense in $B$, we then conclude that $BI \subseteq [IB]$, whence also
  $$
  [BI] \subseteq [IB].
  $$
  The reverse inclusion follows by taking adjoints, so (ii) is proved.  Finally (iii) follows immediately from (i) and (ii).
\endProof

Let us now discuss the naturally arising question as to whether the converse of \ref {JAIsBInvar.ii} holds.
Unfortunately we have no direct proof for this fact but we are nevertheless able to settle it for the special case of
regular ideals under the ideal intersection property:

\state Proposition \label SameInvar
  Let $B$ be a C*-algebra and let $A$ be a regular subalgebra of $B$ satisfying the ideal intersection
property.  Then every $B$-invariant regular ideal of $A$ is normalizer-invariant.

\Proof
Given a $B$-invariant regular ideal  $I\trianglelefteq A$, we have by \ref {MainThm} that
  $$
  I=\alpha (\beta (I)) =   \beta (I)\cap A.
  $$
  Since $\beta (I)$ is an ideal of $B$, the result follows from \ref {JAIsBInvar.i}.
\endProof

  As a consequence of \ref {JAIsBInvar.ii} and \ref {SameInvar} we observe that, under the hypotheses above ($A$ is
regular and satisfies the ideal intersection property), the family of
all regular, $B$-invariant ideals of $A$, denoted by $\Rinv (A)$ in \ref {MainThm}, coincides with the family
  \bitem $\Rprimeinv (A)$ formed by all regular, normalizer-invariant ideals of $A$.
  \hfill \advseqnumbering \label NewA \ref {NewA}

\medskip
We are not aware of any nontrivial example of a subalgebra satisfying axiom \axi /, except for those satisfying the
stronger hypotheses of \ref {JAIsBInvar}.  In fact our motivation has always been the study of regular
subalgebras where axiom \axi / emerged somewhat naturally.

Since axiom \axi / is more general and concise if compared to
regularity, and since it is enough to prove Theorem \ref {MainThm}, we believe its use is partially justified.  A full
justification will of course come with the discovery of a nontrivial and relevant example of a non-regular subalgebra $A$ satisfying axiom \axi /.

The statement of the following immediate consequence of \ref {MainThm}, \ref {JAIsBInvar.iii}, and \ref {SameInvar}
  makes no reference to the perhaps less popular concepts of $B$-invariance and axiom \axi /.
A first version of it appears in \cite [Theorem 3.18]{RegIdeals} under somewhat more stringent  hypotheses.

\state Corollary \label MainCorol
  Let $B$ be a C*-algebra and let $A$ be a regular subalgebra of $B$ satisfying the ideal intersection property.
  Then the conclusion of \ref {MainThm}  holds with $\Rinv (A)$ replaced by $\Rprimeinv (A)$,  as
defined in \ref {NewA}.

\endsection

\startsection The importance of regularity for ideals

\sectiontitle

As we already mentioned in the introduction, the fact that \ref {RenaultCorresondence} is a bijection requires stronger
assumptions as compared to the ideal intersection property, at least in the case of groupoid C*-algebras.
Theorem \ref {MainThm} is therefore not expected to be extendable to
non-regular ideals.  In what follows we describe an elementary example to show that this is indeed the case.

Let $X$ be the topological space $[0,2]$, and let $Y=[0,1]\cup [2,3]$.  Consider the continuous map
  $$
  h: Y\to X
  $$
  given by $h(y)=y$, when $y$ is in $[0,1]$, while $h(y)=y-1$, otherwise.  Thus, $h$ maps the first connected component
of $Y$ bijectively onto the first half of $X$, and it maps the second component of $Y$ bijectively onto the second half
of $X$.  We then see that $h$ is onto $X$, but it is not one-to-one due to the fact that $h(1)=h(2)$.

Let $B=C(Y)$, and let $A=C(X)$,  viewed as a subalgebra of $B$ under the monomorphism
  $$
  f\in A \hookrightarrow f\circ h\in B.
  $$

Since $B$ is commutative and hence all subsets of $B$  whatsoever
are $B$-invariant, we see that  $A$ satisfies axiom \axi /.
Also it
follows from \cite [2.11.14]{ExelPitts} that $A$  satisfies the ideal intersection property.
  Theorem \ref {MainThm}  then applies but we claim that it cannot be extended beyond regular ideals.  To see this we will exhibit
two distinct ideals $J_1$ and $J_2$ of $B$, not all regular, such that $J_1\cap A=J_2\cap A$.

Choosing
  $$
  J_1=\{f\in C(Y): f(1)=0\},
  $$
  and
  $$
  J_2=\{f\in C(Y): f(1)=f(2)=0\},
  $$
  it is easy to see that
  $$
  J_1\cap A=J_2\cap A=\{g\in C(X): g(1)=0\},
  $$
  so the map $\alpha $ referred to in \ref {MainThm} fails to be injective.

\medskip
Returning to the analysis of regular ideals,
observe that, in the present example, there is no conditional expectation from $C(Y)$ onto $C(X)$, so \cite [Theorem
3.18]{RegIdeals} does not apply.   However, as already mentioned, the hypotheses of \ref {MainThm} are satisfied, so
one concludes that the boolean algebra of regular ideals of $C(Y)$ is isomorphic to the boolean algebra of all regular
ideals of $C(X)$ (notice that, due to commutativity, every ideal of $C(X)$ is  $C(Y)$-invariant).

\endsection
\references
\overfullrule 0pt

\Article BCFS
  J. Brown, L. Clark, C. Farthing, and A. Sims;
  Simplicity of algebras associated to \'etale groupoids;
  Semigroup Forum, 88 (2014), no. 2, 433-452

\Bibitem RegIdeals
  J. Brown, A. Fuller, D. Pitts and S. Reznikoff;
  Regular ideals, ideal intersections, and quotients;
  arXiv:2208.09943 [math.OA], 2022

\Article Lots
  J. Brown, G. Nagy, S. Reznikoff, A. Sims and D. Williams;
  Cartan subalgebras in C*-algebras of Hausdorff \'etale groupoids;
  Integral Equations Operator Theory, 85 (2016), no. 1, 109-126

\Article CK
  J. Cuntz and W. Krieger;
  A class of C*-algebras and topological Markov chains;
  Inventiones Math., 56 (1980), 251-268

\Article Dixmier
  J. Dixmier;
  Sous-anneaux abeliens maximaux dans les facteurs de type fini;
  Ann. of Math., 59 (1954), 279-286

\Article ExelNHausd
  R. Exel;
  Non-Hausdorff groupoids;
  Proc. Amer. Math. Soc., 139 (2011), 897-907

\Bibitem ExelPitts
  R. Exel and D. Pitts;
  Characterizing groupoid C*-algebras of non-Hausdorff \'etale groupoids;
  Lecture Notes in Mathematics, 2306. Springer, Cham, 2022. viii+158 pp

\Bibitem Renault
  J. Renault;
  A groupoid approach to C*-algebras;
  Lecture Notes in Mathematics vol.~793, Springer, 1980

\Article RenaultCartan
  J. Renault;
  Cartan subalgebras in C*-algebras;
  Irish Math. Soc. Bulletin, 61 (2008), 29-63

\endgroup

\close